\documentclass{amsart}

\usepackage{graphicx}

\usepackage{lineno}

\usepackage{hyperref}
\hypersetup{hidelinks=true}

\newtheorem{theorem}{Theorem}[section]
\newtheorem{lemma}[theorem]{Lemma}

\theoremstyle{definition}
\newtheorem{definition}[theorem]{Definition}
\newtheorem{corollary}[theorem]{Corollary}

\newtheorem{proposition}[theorem]{Proposition}

\theoremstyle{remark}

\numberwithin{equation}{section}

\newcommand{\al}{\alpha}
\newcommand{\bb}{\beta}

\newcommand{\Z}{\mathbb Z}

\newcommand{\T}{\mathcal{T}}
\newcommand{\R}{\mathcal{R}}
\newcommand{\s}{\mathcal{S}}
 
 \newcommand{\ip}[1]{\left\langle#1\right\rangle}

\usepackage[dvipsnames]{xcolor}
\definecolor{matblue}{rgb}{0, 0.447, 0.741}

\begin{document}

\title{A tree distinguishing polynomial}


\author{}
\address{}
\curraddr{}
\email{}
\thanks{Department of Mathematics, Simon Fraser University. Burnaby, BC V5A 1S6, Canada.}
\thanks{E-mail: \href{mailto:pengyu_liu@sfu.ca}{pengyu\_liu@sfu.ca}}

\author{Pengyu Liu}
\address{}
\curraddr{}
\email{}
\thanks{}



\date{}

\dedicatory{}

\begin{abstract}
We define a bivariate polynomial for unlabeled rooted trees and show that the polynomial of an unlabeled rooted tree $T$ is the generating function of a class of subtrees of $T$. We prove that the polynomial is a complete isomorphism invariant for unlabeled rooted trees. Then, we generalize the polynomial to unlabeled unrooted trees and we show that the generalized polynomial is a complete isomorphism invariant for unlabeled unrooted trees. 
\end{abstract}

\maketitle


\section{Introduction}

\noindent
Polynomial invariants are important tools in the study of graphs, knots and links. The Tutte polynomial \cite{T} is the most investigated polynomial invariant for graphs. As an isomorphism invariant of graphs, the Tutte polynomial carries some information of a graph, for example, the chromatic polynomial and the number of spanning trees of the graph. The well known Jones polynomial \cite{Jones} and the HOMFLY polynomial \cite{HOMFLY} are important invariants for knots and links which are related to the crossing number and the braid index of knots and links respectively \cite{Adams}. However, the Tutte polynomial fails to distinguish trees. Actually, all trees with the same number of edges have the same Tutte polynomial. In the doctoral thesis of Law \cite{Law}, polynomials were divided into three levels according to their tree distinguishing power, where the most powerful (level three) polynomial is the chromatic symmetric function introduced in 1995 by Stanley \cite{Stan}.  It is proved that the chromatic symmetric function distinguishes some classes of trees including spiders \cite{Mar} and caterpillars \cite{Loebl}, but it remains a conjecture that the chromatic symmetric function is a complete isomorphism invariant for trees. The $U$-polynomial defined by Noble and Welsh \cite{NW}, which is equivalent to the polychromate introduced by Brylawski \cite{Bry,MN,Sa}, determines the chromatic symmetric function and vice versa when restricted to trees. See \cite{NW} and \cite{A}. The strong polychromate defined by Bollob\'as and Riordan \cite{Bo} is also equivalent to the polynomials above when restricted to trees. Hence, whether these polynomials are complete invariants for trees depends on whether the chromatic symmetric function is a complete invariant for trees. The level two polynomials are the polynomials that distinguish rooted trees. These polynomials include the subtree polynomial introduced by Chaudhary and Gordon \cite{Chow}, the Ising polynomial introduced by  Andr\'en and Markstr\"om \cite{Ising} and the Negami polynomial introduced by Negami and Ota \cite{Negami}. However, it has been unknown to date whether there exists a polynomial that is a complete isomorphism invariant for unrooted trees. 

\bigskip
\noindent
In the emerging fields of phylogenetics and linguistics, the information carried by the shapes of trees needs to be analyzed and compared quantitively and accurately. A polynomial invariant, especially a complete invariant, for trees is a potentially convenient tool for this task because polynomials are well studied mathematical objects. With this motivation, we introduce a new polynomial that is a complete isomorphism invariant for trees, which, to our knowledge, is the first of its kind.

\bigskip
\noindent 
Before demonstrating the results, we clarify the terminology used in this paper. Trees are unlabeled unless otherwise stated. The order of a tree stands for the number of vertices in the tree.  A monomial always has coefficient one and a term in a polynomial may have any integer coefficient. All of the internal vertices of a rooted $m$-ary tree have exactly $m$ children. Similarly, every internal vertex of an unrooted $m$-ary tree has degree $m+1$. This paper is structured as follows. First, we define a polynomial for rooted trees and investigate the information about the trees carried by the polynomial. Then, we show that the polynomials for rooted trees are irreducibles in the polynomial ring and that the polynomial distinguishes rooted trees. Finally, we introduce a systematic way to generalize a polynomial that distinguishes rooted trees to a polynomial that distinguishes unrooted trees.

\section{A polynomial for rooted trees}

\subsection{Definitions}\label{r1}

Let $\T_r$ be the set of rooted trees and $k\geq 1$ be an integer. The {\em rooted $k$-star} is the tree in $\T_r$ of order $k+1$ with $k$ leaf vertices in which all of the leaf vertices are adjacent to the single internal vertex that is identified as the root. The {\em $k$-wedge operation} $\wedge_k:\T_r^k\to\T_r$ is defined such that for $k$ rooted trees $T_1$, $T_2$,...,$T_k$, $\wedge_k(T_1,T_2,...,T_k)$ is the tree in $\T_r$ constructed by pasting the roots of the $k$ trees to the $k$ leaf vertices of the rooted $k$-star respectively. Note that the operation $\wedge_k$ is {\em permutative}, that is, for any permutation $\pi\in S_k$, where $S_k$ is the symmetric group, $\wedge_k(T_1,T_2,...,T_k)=\wedge_k(T_{\pi(1)},T_{\pi(2)},...,T_{\pi(k)})$. Besides, the rooted $m$-ary trees can be constructed recursively using only $\wedge_m$. In particular, the rooted binary trees can be constructed recursively using only $\wedge_2$.

\begin{definition}\label{df}
Let $\bullet$ be the trivial tree with one vertex and $T=\wedge_k(T_1,T_2,...,T_k)$ be an arbitrary tree in $\T_r$, where $k\geq 1$. The polynomial $P:\T_r\to \Z[x,y]$ is defined by the following rules.
\begin{enumerate}
\item $P(\bullet)=x$,
\item $P(T)=y+\prod_{i=1}^k P(T_i)$.
\end{enumerate}
\end{definition}

\noindent
For example, the rooted $k$-star $T$ can be constructed by applying $\wedge_k$ to $k$ trivial trees, so its polynomial is $P(T)=y+x^k$. Let $l\geq0$ be an integer. The {\em rooted path} of length $l$ is the path with $l+1$ vertices such that one of the leaf vertices is identified as the root. The rooted path $T$ of length $l$ can be constructed by applying  the $1$-wedge operation $l$ times starting from the trivial tree. According to the definition above, the polynomial of $T$ is $P(T)=ly+x$. 

\bigskip
\noindent
Let $T$ be a rooted tree in $\T_r$ and $v$ be an vertex of $T$. The {\em affix tree} of $T$ to the vertex $v$, denoted by $T'_v$, is the subtree of $T$ induced from the vertex $v$ and all of the descendants of $v$. Note that if $v$ is a leaf vertex, then $T'_v$ is the trivial tree with a single vertex. If the root vertex of $T$ has degree one, the {\em branching vertex} of $T$ is defined to be the nearest vertex to the root with degree greater than two. If $T$ does not have a vertex with degree greater than two, then it is a rooted path and it does not have a branching vertex. If the root vertex of $T$ has degree greater than one, then we consider that the root vertex is also the branching vertex. In this paper, the affix tree of $T$ to the branching vertex will be mentioned many times, hence, we denote this specific affix tree by $T'$. The {\em stem} of $T$ is the path from its root vertex to its branching vertex. See Figure \ref{f1} for an example. In particular, if the root vertex of $T$ has degree greater than one, that is, the root vertex is also the branching vertex, then the stem consists of only the root vertex and is of length zero.

\begin{proposition}\label{welldefine}
The function $P:\T_r\to \Z[x,y]$ is well defined.
\end{proposition}
\proof Denote the number of leaf vertices of a tree in $\T_r$ by $n$. We prove the proposition by strong induction on $n$.
\begin{enumerate}
\item If $n=1$,  the rooted trees in $\T_r$ with a single leaf vertex are rooted paths. The polynomial of the rooted path $T$ of length $l$ is $P(T)=ly+x$. If two rooted paths are isomorphic, they must have the same length, hence the same polynomial.
\item Assume that for $n\leq N$, the polynomials of isomorphic trees are identical.
\item If $n=N+1$, let $T$ and $B$ be two isomorphic trees in $\T_r$. $T \simeq B$ implies that their stems are of the same length and their affix trees to the branching vertices $T'=\wedge_k(T_1, T_2,...,T_k)$ and  $B'=\wedge_l(B_1,B_2,...,B_l)$ are also isomorphic, hence $k=l>1$. Note that the numbers of leaf vertices in $T_i$ and $B_i$ are fewer than $N+1$ for all $1\leq i\leq k$. Without loss of generality, we will compare $T_i$ with $B_i$ since the $k$-wedge operations are permutative. $T'\simeq B'$ implies $T_i\simeq B_i$ for all $1\leq i\leq k$. According to the hypothesis, $P(T')=y+\prod_{i=1}^k P(T_i)=y+\prod_{i=1}^k P(B_i)=P(B')$. Assuming that the stems of $T$ and $B$ are of length $l$, we can construct $T$ and $B$ by recursively applying the $1$-wedge operation $l$ times to $T'$ and $B'$.  According to Definition \ref{df}, $P(T)=ly+P(T')=ly+P(B')=P(B)$. \qed 
\end{enumerate} 

\noindent
Note that according to the proof above, a tree $T$ in $\T_r$ with a stem of length $l$ has the polynomial $P(T)=ly+P(T')$, where $T'$ is the affix tree of $T$ to the branching vertex.
 
\bigskip
\noindent
To compute the polynomial of a tree in $\T_r$, we can apply the recurrence relation in Definition \ref{df}. Besides, the polynomial can also be computed directly using the Dyck word \cite{Ghys} of the tree by placing a symbol $x$ in every pair of parentheses that represents a leaf vertex and placing the symbol $+y$ before the end of every pair of parentheses that represents an internal vertex. For example, if the Dyck word of a tree is $((()()))$ then its polynomial should be $(((x)(x)+y)+y)=x^2+2y$. On the other hand, to reconstruct a tree from its polynomial, we may compute its Dyck word by recursively subtracting $y$ and factoring the rest of the polynomial. See \cite{Mon} for methods to factor large multivariate polynomials.

\subsection{Interpretation of the polynomial}\label{ioc}
\begin{figure}[!htb]
\begin{center}
\includegraphics[scale=0.5]{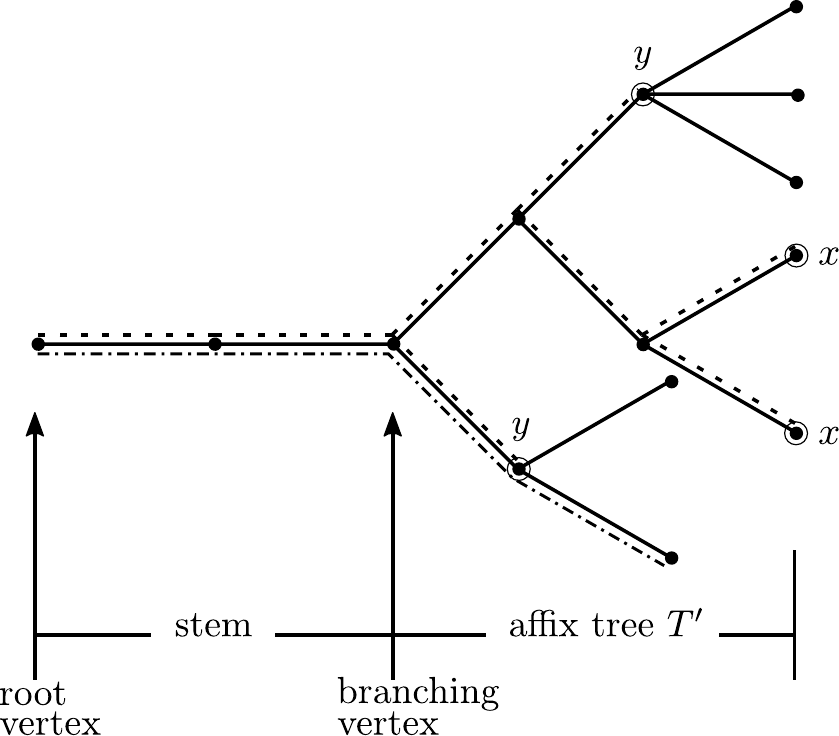}
\caption{An example of a rooted tree $T$ where dashed lines represent a primary subtree $S$ of $T$ whose monomial $q(S)=x^2y^2$ and the dash-dotted line represents a subtree of $T$ that is not a primary subtree.}\label{f1}
\end{center}
\end{figure}

To determine how the polynomial describes the features of trees in $\T_r$, we introduce the following definitions. Let $T$ be a tree in $\T_r$.  A {\em primary subtree} $S$ of $T$ is a rooted subtree of $T$ such that $S$ shares the same root vertex with $T$ and any leaf vertex of $T$ is either a leaf vertex of $S$ or a descendant of a leaf vertex of $S$. In other words, no leaf vertex of $T$ can be a descendant of an internal vertex of a primary subtree $S$ of $T$. Denote the set of all primary subtrees of $T$ by $\s_T$. For any primary subtree $S$ of $T$, we assign a monomial $q(S)=x^\al y^\bb$ to the primary subtree if $S$ possesses $\al+\bb$ leaf vertices and $\al$ of them are leaf vertices of $T$ and $\bb$ of them are internal vertices of $T$. Hence, the total degree of $q(S)$ is the number of leaf vertices of $S$. See Figure \ref{f1} for an example. Note that if a tree in $\T_r$ is not the trivial tree, we always consider the root vertex of the tree as an internal vertex even if it has degree one. For example, the rooted $k$-star has two primary subtrees. One is the trivial tree with only the root vertex, which corresponds to a monomial $q(S)=y$. The other is the tree that is isomorphic to the rooted $k$-star, which corresponds to the monomial $q(S)=x^k$. The rooted path $T$ of length $l$ has $l+1$ primary subtrees, which are the paths from the root to each of its vertices including the root. If $S$ is the primary subtree from the root to the leaf vertex, then $q(S)=x$ since the leaf vertex of $S$ is also the leaf vertex of $T$. For any other primary subtree $S$ of $T$, $q(S)=y$ because the leaf vertex of any $S$ is an internal vertex of $T$. 

\bigskip
\noindent
Let $G_1=(V_1,E_1)$ and $G_2=(V_2,E_2)$ be two subgraphs of a graph $G$ and $V_1$, $V_2$, $E_1$, $E_2$ are the sets of vertices and edges of $G_1$ and $G_2 $ respectively. We define the intersection $G_1\cap G_2=(V_1\cap V_2,E_1\cap E_2)$ and $G_1\cap G_2=\emptyset$ if $V_1\cap V_2=\emptyset$ and $E_1\cap E_2=\emptyset$. 

\begin{lemma}\label{addlemma} \em
Let $T$ be a tree in $\T_r$ and $T'$ be the affix tree of $T$ to the branching vertex. The following statements about primary subtrees and the monomials are true.
\begin{enumerate}
\item For any primary subtree $S$ of $T$, if $S\cap T'\neq\emptyset$, then $S'=S\cap T'$ is a primary subtree of $T'$ and $q(S)=q(S')$.
\item If $T=\wedge_k(T_1,T_2,...,T_k)$ with $k>1$, any primary subtree $S$ of $T$ is of form $\wedge_k(S_1,S_2,...,S_k)$, where $S_i\in\s_{T_i}$ for all $1\leq i\leq k$, except for the primary subtree $S_r$ which consists of only the root vertex of $T$.
\item Suppose that $T=\wedge_k(T_1,T_2,...,T_k)$ with $k>1$ and $S$ be a primary subtree of $T$ such that $S=\wedge_k(S_1,S_2,...,S_k)$, where $S_i\in\s_{T_i}$ and $q(S_i)=x^{\al_i}y^{\bb_i}$ for all $1\leq i\leq k$. Then $q(S)=\prod_{i=1}^k q(S_i)=x^{\sum_{i=1}^k \al_i}y^{\sum_{i=1}^k \bb_i}$ which is still a monomial. 
\end{enumerate}
\end{lemma}

\proof 
For (1), if $S\cap T'=\emptyset$, then $S$ is a primary subtree on the stem of $T$ and $q(S)=y$. This is because $S'=S\cap T'\neq\emptyset$ is a rooted tree that shares the same root with $T'$ and any leaf vertex or internal vertex of $T'$ is also a leaf vertex or an internal vertex of $T$ respectively, hence, any leaf vertex of $T'$ is a leaf vertex of $S'$ or a descendant of a leaf vertex of $S'$ and according to the definition of a primary subtree and the definition of the monomial, $S'$ is a primary subtree of $T'$ with $q(S)=q(S')$. If $S\cap T'=\emptyset$, then $S$ contains no vertex of $T'$, that is, all the vertices of $S$ are on the stem. Therefore, $S$ is a rooted path and its leaf vertex is an internal vertex of $T$ and $q(S)=y$. 

For (2), Note that any primary subtree $S\not\simeq S_r$ of $T$ is rooted at the root vertex of $T$ and every leaf vertex of $T$ is a descendant of a leaf vertex of $S$. Hence, for any $1\leq i\leq k$, there exists at least one leaf vertex of $S$ in $T_i$. Let $L_i$ be the set of leaf vertices of $S$ in $T_i$ and $S_i$ be the induced subtree of $T_i$ from $L_i$ and all the ancestors of vertices in $L_i$. $S_i$ is a primary subtree of $T_i$ since $S$ is a primary subtree of $T$ and all the leaf vertices of $T_i$ are descendants of vertices in $L_i$. Therefore, the primary subtree $S=\wedge_k(S_1,S_2,...,S_k)$. Conversely, if $T=\wedge_k(T_1,T_2,...,T_k)$ with $k>1$, any tree $S$ of form $\wedge_k(S_1,S_2,...,S_k)$, where $S_i\in\s_{T_i}$ for all $1\leq i\leq k$, is a primary subtree of $T$ according to the definition of primary subtrees. 

For (3), the leaf vertices of $S_i$ that are leaf vertices of $T_i$ are also leaf vertices of $T$ and the leaf vertices of $S_i$ that are internal vertices of $T_i$ are also internal vertices of $T$ for any $1\leq i\leq k$. \qed

\begin{lemma}\label{lemma} 
$P(T)=\sum_{S\in\s_T}q(S)$.
\end{lemma}

\proof 
We prove the lemma by strong induction on $n$, the number of leaf vertices of a tree in $\T_r$.
\begin{enumerate}
\item If $n=1$, we know, from Section \ref{r1}, that the rooted trees in $\T_r$ with a single leaf vertex are rooted paths and the polynomial of the rooted path $T$ of length $l$ is $P(T)=ly+x$. On the other hand, we also know that $T$ has $l+1$ primary subtrees. The primary subtree $S$ that is isomorphic to $T$ has the monomial $q(S)=x$ and any of the other $l$ primary subtrees has the monomial $q(S)=y$. Hence, $P(T)=\sum_{S\in\s_T}q(S)$.
\item Assume $P(T)=\sum_{S\in\s_T}q(S)$ for all trees in $\T_r$ with $n\leq N$ leaf vertices.
\item If $n=N+1$, let $T$ be an arbitrary tree in $\T_r$ with $N+1$ leaf vertices. Suppose $T$ has a stem of length $l$ and its affix tree to the branching vertex $T'=\wedge_k(T_1,T_2,...,T_k)$, where $k>1$ and $T_i$ is a tree in $\T_r$ with fewer than $N+1$ leaf vertices for all $1\leq i\leq k$. We know that $P(T)=ly+P(T')$ according to Section \ref{r1}. Note that the first fact above states that there is a partition of $\s_T$ such that for any $S\in\s_T$, if $S\cap T'=\emptyset$, then $S$ is a rooted path with vertices on the stem of $T$. There exist $l$ such primary subtrees of $T$, namely, the paths from the root to each of the vertices of the stem except for the branching vertex, which contributes to the $ly$ term in the polynomial $P(T)$. Besides, for any $S\in\s_T$, if $S'=S\cap T'\neq\emptyset$, then $S'\in\s_{T'}$ and $q(S)=q(S')$. Therefore, if we can prove that $P(T')=\sum_{S\in\s_{T'}}q(S)$, then $P(T)=\sum_{S\in\s_T}q(S)$ follows. According to Definition \ref{df} and the induction hypothesis, $P(T')=y+\prod_{i=1}^k P(T_i)=y+\prod_{i=1}^k\sum_{S\in\s_{T_i}}q(S)$. The monomial $y$ in $P(T')$ corresponds to the primary subtree $S_r$ of $T'$ since it is the only primary subtree of $T$ that has one leaf vertex in $T'$. For any monomial in $\prod_{i=1}^k\sum_{S\in\s_{T_i}}q(S)$, it is of form $q(S_1)q(S_2)...q(S_k)$ where $S_i\in\s_{T_i}$ for any $1\leq i\leq k$. We know from the second and the third facts above that $q(\wedge_k(S_1,S_2,...,S_k))=\prod_{i=1}^k q(S_i)$ and $\wedge_k(S_1,S_2,...,S_k)$ is a primary subtree of $T'$. So every monomial in $P(T')$ is a monomial in $\sum_{S\in\s_{T'}}q(S)$. On the other hand, for any $S\in\s_{T'}$ except for the trivial primary subtree $S_r$, $S=\wedge_k(S_1,S_2,...,S_k)$ where $S_i\in\s_{T_i}$ for any $1\leq i\leq k$ according to the second and the third facts above. For the primary subtree $S_r\in\s_{T'}$, $q(S_r)=y$ and $y$ is a monomial in $P(T')$. For any other primary subtree $S\in\s_{T'}$, $q(S)=\prod_{i=1}^k q(S_i)$ which is also a monomial in $P(T')$. Therefore, $P(T')=\sum_{S\in\s_{T'}}q(S)$ and $P(T)=\sum_{S\in\s_T}q(S)$. \qed 
\end{enumerate}

\bigskip
\noindent 
Lemma \ref{lemma} shows that the polynomial of a tree $T$ can be interpreted as the generating function of the number of primary subtrees whose set of leaf vertices consists of $\al$ leaf vertices of $T$ and $\bb$ internal vertices of $T$. 
\begin{corollary}
Let $T$ be a tree in $\T_r$. Suppose $P(T)$ has $m$ terms and $a_1,a_2,...,a_m$ are the corresponding coefficients, then $T$ has $\sum_{i=1}^{m} a_i$ primary subtrees.
\end{corollary}

\noindent
Let $T$ be a rooted tree with $n$ leaf vertices in $\T_r$. According to Lemma \ref{lemma}, there exists a $x^n$ term in $P(T)$ which corresponds to the primary subtree that is isomorphic to $T$, that is, all the leaf vertices of $T$ are leaf vertices of $S$. There exists only one such primary subtree, hence the coefficient of the term is one. Moreover, for any other primary subtree $S$ of $T$, $q(S)$ has a factor $y$ because at least one leaf vertex of $T$ is not a leaf vertex of $S$, that is, at least one leaf vertex of $T$ is a descendant of a leaf vertex of $S$ which is an internal vertex of $T$. Last but not least, there exists at least one primary subtree with only one leaf vertex. The subtree $S_r$ of $T$ consisting of only the root vertex is such a primary subtree and it exists for any tree in $\T_r$. Note that the leaf vertex of such primary subtrees is always an internal vertex of $T$ except for the trivial tree. Therefore, if $T$ is not the trivial tree, then there is always a term $ty$ in the polynomial $P(T)$, where $t$ is the number of primary subtrees with only one leaf vertex which is an internal vertex of $T$. Indeed, $t=l$ if $T$ is a rooted path and $t=l+1$ otherwise, where $l$ is the length of the stem of $T$. Besides, no primary subtrees of $T$ other than those on the stem can contribute to a monomial $y$ since if a primary subtree $S$ of $T$ contains vertices that are not on the stem of $T$, $S$ must have at least two leaf vertices.

\subsection{Complete isomorphism invariants for rooted trees}\label{r3}
To prove that the polynomial $P:\T_r\to \Z[x,y]$ is a complete isomorphism invariant for rooted trees, we need to prove the polynomials are irreducibles in $\Z[x,y]$. Eisenstein's criterion states that if $D$ is an integral domain, $I_p$ is a prime ideal of $D$ and $Q=a_nx^n+a_{n-1}x^{n-1}+...+a_1x+a_0$ is a polynomial in $D[x]$, then $Q$ is an irreducible in $D[x]$ if $a_n\not\in I_p$, for all $0\leq i<n$ $a_i\in I_p$ and $a_0\not\in I_p^2$.

\begin{lemma}\label{irredu} \em
For any tree $T$ in $\T_r$, $P(T)$ is an irreducible in $\Z[x,y]$.
\end{lemma}
\proof We use Eisenstein's criterion to prove this lemma. Note that $\Z[x,y]=\Z[y][x]$. Let $\Z[y]$ be the integral domain $D$ and $I_p=\ip{y}$ be the prime ideal in $\Z[y]$. Suppose $T$ has $n$ leaf vertices. We know, from Section \ref{ioc}, that the leading term of $P(T)$ is always $x^n$ and $a_n=1$, hence, $a_n\not\in I_p$. For any primary subtree $S$ that is not isomorphic to $T$, there is always a leaf vertex of $S$ that is an internal vertex of $T$, so $q(S)$ always has a factor $y$ and $a_i\in I_p$ for all $0\leq i<n$. Moreover, the constant term $a_0$ always contains a term $ty$. Therefore, $a_0\not\in I_p^2$ and the polynomials  for all rooted trees in $\T_r$ are irreducibles in $\Z[x,y]$. \qed

\begin{proposition}\label{1-1}
The function $P:\T_r\to \Z[x,y]$ is injective.
\end{proposition}
\proof We prove the proposition by strong induction on $n$, the number of leaf vertices of a tree in $\T_r$.
\begin{enumerate}
\item  If $n=1$, the polynomial of a rooted path $T$ of length $l$ is $P(T)=ly+x$. Two non-isomorphic rooted paths have different lengths, so their polynomials are different.
\item Assume the function is injective for all $n\leq N$.
\item If $n=N+1$, let $T$ and $B$ be two non-isomorphic trees in $\T_r$ with $N+1$ leaf vertices. Note that only the primary subtrees on the stem of a tree in $\T_r$ contribute to the $ty$ term. If the stems of $T$ and $B$ are of different lengths, the $ty$ terms in the polynomials of $T$ and $B$ will have different coefficients hence $P(T)\neq P(B)$. Suppose that the stems of $T$ and $B$ are of the same length, $T'=\wedge_{k}(T_1,T_2,...,T_k)$ and $B'=\wedge_l(B_1,B_2,..,B_l)$ where for all $1\leq i \leq k$, $T_i$ and $B_i$ are rooted trees in $\T_r$ with fewer than $N+1$ leaf vertices. If $P(T)=P(B)$, then $P(T')=y+\prod_{i=1}^kP(T_i)=y+\prod_{j=1}^lP(B_i)=P(B')$, that is, $P(T_1)P(T_2)...P(T_k)=P(B_1)P(B_2)...P(B_l)$.  Since these polynomials are irreducibles in $\Z[x,y]$ according to Lemma \ref{irredu} and $\Z[x,y]$ is a unique factorization domain, we know $k=l$ and, without loss of generality,  $P(T_i)=P(B_i)$ for all $1\leq i\leq k$ after a rearrangement of labels. Then, the hypothesis implies that $T_i\simeq B_i$ for all $1\leq i\leq k$, hence, $T'\simeq B'$. Note that the stems of $T$ and $B$ are of the same length. Therefore, $T\simeq B$ which contradicts the assumption. \qed

\end{enumerate}

\bigskip
\noindent 
Proposition \ref{welldefine} and Proposition \ref{1-1} imply that the polynomial is a complete isomorphism invariant for rooted trees.

\begin{theorem}\label{main}
$T_1, T_2\in\T_r$ are isomorphic if and only if $P(T_1)=P(T_2)$. 
\end{theorem}

\noindent
Let $\T_r^*$ be the set of rooted trees such that every internal vertex has more than one child and $\T_r^m$ be the set of rooted $m$-ary trees where $m\geq2$ is an integer. For any tree $T$ in $\T_r$ and any prime number $p$, the polynomial $P_p:\T_r\to\Z[x]$ is defined by substituting  $p$ for $y$ in the polynomial $P(T)$. We can prove that the polynomial  $P_p(T)$ is an irreducible in $\Z[x]$ for any tree $T$ in $\T_r$ by substituting a prime number $p$ for $y$ in the proof of Lemma \ref{irredu}. Then, the following corollary follows the proof of Proposition~\ref{1-1}.

\begin{corollary}\label{prime}
$T_1, T_2\in\T_r^*$ are isomorphic if and only if $P_p(T_1)=P_p(T_2)$. In particular,  $T_1, T_2\in\T_r^m$ are isomorphic if and only if $P_p(T_1)=P_p(T_2)$
\end{corollary}

\begin{figure}[!htb]
\begin{center}
\includegraphics[scale=0.5]{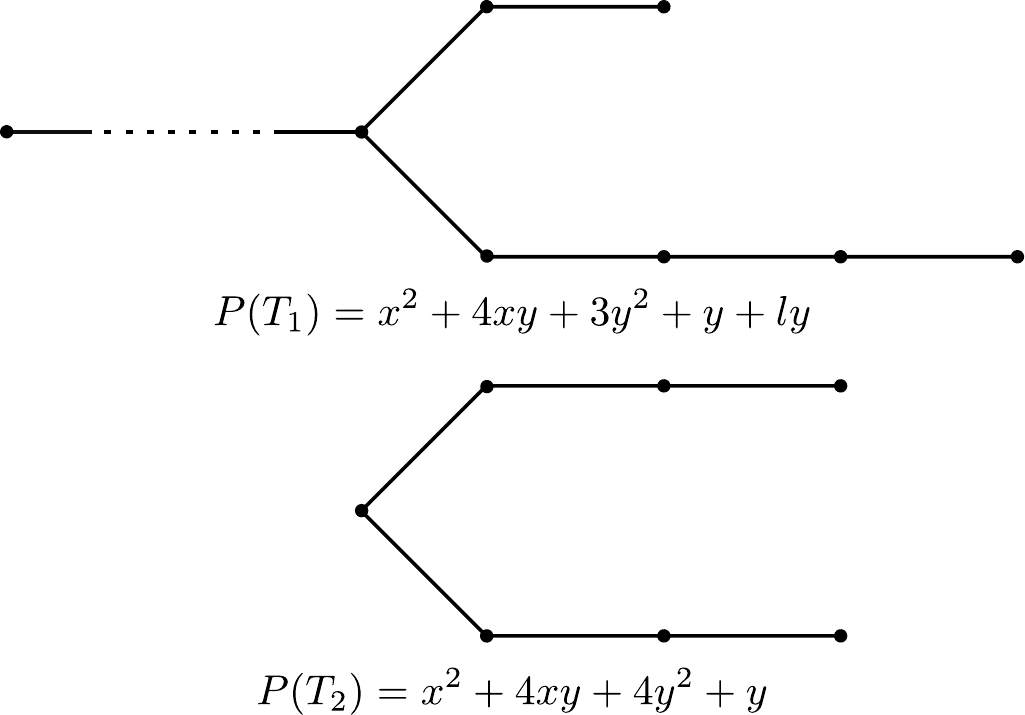}
\caption{A pair of rooted trees that share the same polynomial $P_p(T):\T_r\to\Z[x]$.}\label{f2}
\end{center}
\end{figure}

\noindent However, for any prime number $p$, there exists a pair of non-isomorphic rooted trees in $\T_r$ with the same polynomial $P_p(T):\T_r\to\Z[x]$, where there exists at least one internal vertex that has only one child. Figure \ref{f2} shows a pair of such trees. For any prime number $p$, we can choose the length of the stem of the tree $T_1$ to be $l=p$  $P_p(T_1)=P_p(T_2)$.  An interesting question is determine the values of the integer $n$ such that $P_n(T):\T_r\to\Z[x]$ is a complete isomorphism invariant for rooted $m$-ary trees or trees in $\T_r^*$. 
If $m=2$, it can be checked by computer that $P_n(T):\T_r\to\Z[x]$ is not a complete isomorphism invariant for rooted $2$-ary or binary trees when $n\in\{-1,0,1\}$. For any other integer, it is not known whether $P_n(T):\T_r\to\Z[x]$ is a complete isomorphism invariant for rooted binary trees or not. It is not known either for other rooted $m$-ary trees. 

\section{A tree distinguishing polynomial}

\subsection{The polynomial for unrooted trees}\label{r4}
Let $\T_u$ be the set of unrooted trees, and $\T=\T_u\cup\T_r$. Suppose $T$ is a tree in $\T_u$ with $n$ leaf vertices. A {\em leaf edge} of $T$ is an edge of $T$ that is incident to a leaf vertex. For each leaf edge of $T$, we can construct a rooted tree $T_i$ by contracting the leaf edge and identifying the contracted edge as the root vertex of $T_i$. Denote the set of such rooted trees constructed from $T$ by $\R_T$. Note that $\R_T$ has $n$ elements and some of them may be isomorphic.

\begin{lemma}\label{ll}\em
$T, B\in\T_u$ are isomorphic if and only if there exists a bijection $h:\R_{T}\to\R_{B}$ such that for any $T_i$ in $\R_T$, $T_i$ is isomorphic to $h(T_i)$.
\end{lemma}
\proof Suppose $T\simeq B$ and $\phi:T\to B$ is the isomorphism. Let $T_i$ be an arbitrary tree in $\R_T$ and $e$ be the edge of $T$ that is contracted to attain $T_i$. We define a function $h:\R_{T}\to\R_{B}$ such that $h(T_i)=B_j$ where $B_j$ is the tree in $\R_B$ constructed by contracting $\phi(e)$. $h:\R_{T}\to\R_{B}$ is a bijection because $\phi:T\to B$ is a bijection between the set of leaf edges of $T$ and the set of leaf edges of $B$. Conversely, Suppose that $T_i$ is a rooted tree in $\R_{T}$ and $B_j=h(T_i)$ is a rooted tree in $\R_{B}$. We can reconstruct $T$ and $B$ from $T_i$ and $B_j$ by recovering the contracted edges, that is, adding an edge and a leaf vertex to the root vertices of $T_i$ and $B_j$ respectively. Therefore, $T_i\simeq B_j$ implies $T\simeq B$. \qed

\bigskip
\noindent
Now, we generalize the polynomial in Definition \ref{df} to $P:\T\to\Z[x,y]$ in the following way. If a tree $T$ in $\T$ is rooted, then $P(T)$ is the polynomial defined as in Definition \ref{df}. If a tree $T$ is unrooted, then we define $P(T)=\prod_{T_i\in\R_T}P(T_i)$. For example, the polynomial of the unrooted $3$-star is $(x^2+y)^3$. Note that for any trees $T_1$ and $T_2$ in $\T$, if $T_1$ is rooted and $T_2$ is unrooted, we always consider that $T_1$ is not isomorphic to $T_2$ even if the only difference between $T_1$ and $T_2$ is an identified rooted vertex. We prove that the polynomial  $P: \T\to\Z[x,y]$ is a complete isomorphism invariant for trees. 

\begin{theorem}\label{t2}\em
$T_1, T_2\in\T$ are isomorphic if and only if $P(T_1)=P(T_2)$. 
\end{theorem}
\proof If $T_1\simeq T_2$, then either both of them are rooted or both of them are unrooted. If both of them are rooted, then $P(T_1)=P(T_2)$ follows Theorem \ref{main}. If both of them are unrooted, it follows from Lemma \ref{ll} and Theorem \ref{main} that $P(T_1)=P(T_2)$. On the other hand, if $T_1\not\simeq T_2$, we have three cases. First, if both of them are rooted, then $P(T_1)\neq P(T_2)$ follows from Theorem \ref{main}. Second, if one of them is rooted and the other is unrooted, then $P(T_1)\neq P(T_2)$ because the polynomial for the rooted tree is an irreducible in $\Z[x,y]$ and the polynomial for the unrooted one is not. Third, if both of them are unrooted, then $P(T_1)\neq P(T_2)$ because otherwise Theorem \ref{main} and $\Z[x,y]$ being a unique factorization domain imply that there exists a bijection $h:\R_{T_1}\to\R_{T_2}$ such that $T\simeq h(T)$ for any $T\in \R_{T_i}$. This contradicts Lemma \ref{ll}, hence, $P(T_1)\neq P(T_2)$. \qed

\bigskip
\noindent
The proof of Theorem \ref{t2} shows that whenever we have a polynomial that represents a class of rooted trees, if (i) the polynomial ring is a unique factorization domain, (ii) the polynomial is a complete isomorphism invariant for the class of rooted trees and (iii) the polynomials of rooted trees in the class are irreducibles in the polynomial ring, then we can generalize the polynomial to the corresponding class of unrooted trees and the resulting polynomial distinguishes these unrooted trees. In particular,  the univariate polynomial $P_p(T):\T_r^m\to\Z[x]$ for rooted $m$-ary trees can be generalized to distinguish $m$-ary trees. Let $\T^m$ be the set of $m$-ary trees including the rooted trees and the unrooted trees and the polynomial $P_p(T):\T^m\to\Z[x]$ is defined such that for any $T$ in $\T^m$, if $T$ is rooted, then $P_p(T)$ is defined as in Section \ref{r3} and if $T$ is unrooted, then $P_p(T)=\prod_{T_i\in\R_T}P_p(T_i)$.

\begin{corollary}
$T_1, T_2\in\T^m$ are isomorphic if and only if $P_p(T_1)=P_p(T_2)$. 
\end{corollary}

\noindent
Now we know that one variable is sufficient to uniquely represent $m$-ary trees by polynomials. An interesting question is whether one variable is sufficient to uniquely represent all trees by polynomials. The zero loci of the polynomials $P: \T\to\Z[x,y]$ for trees may also be interesting.

\subsection{A generalization}
A polynomial distinguishing leaf labeled trees has various applications in linguistics and mathematical biology especially in phylogenetics. The coefficients of a polynomial can be considered as a vector, so norms hence metrics of trees can be induced from tree distinguishing polynomials. Tree metrics, especially metrics for leaf labeled tree, have several biological applications, for example, to compare and classify phylogenetic tree shapes \cite{Colijn}. The polynomial $P:\T_r\to \Z[x,y]$ can be generalized to represent leaf labeled rooted trees in a natural way.  Given a tree $T\in\T_r$ and its polynomial $P(T)$, we can consider the tree $T$ as a vertex labeled tree such that each of its leaf vertices has a label $x$ and an internal vertex $v$ has the polynomial $P(T'_v)$ as its label, where $T'_v$ is the affix tree of $T$ to the vertex $v$. Thus the root vertex of $T$ has the label $P(T)$. If the leaf vertices of a tree have different labels and $\T_r^{\ell}$ denotes the set of leaf labeled rooted trees, we define an analogous polynomial as follows.

\begin{definition}\label{polyll}
Let $\bullet_i$ be the trivial tree with a single vertex that is labeled by $i$ and $T=\wedge_k(T_1,T_2,...,T_k)$ be an arbitrary tree in $\T_r^{\ell}$, where $k\geq 1$. The polynomial $P_{\ell}:\T_r^{\ell}\to \Z[x_1,x_2,...,x_t,y]$ is defined by the following rules.
\begin{enumerate}
\item $P_{\ell}(\bullet_i)=x_i$,
\item $P_{\ell}(T)=y+\prod_{i=1}^k P_{\ell}(T_i)$
\end{enumerate}
\end{definition}

\noindent
Note that different leaf vertices may have the same label. If we set $x_i=x$ for all $1\leq i\leq t$, then $P_{\ell}(T)=P(T)$. The polynomial $P_{\ell}:\T_r^{\ell}\to \Z[x_1,x_2,...,x_t,y]$ is a complete isomorphism invariant for leaf labeled rooted trees, where two leaf labeled trees in $\T_r^{\ell}$ being isomorphic means not only that the unlabeled trees are isomorphic but also that the labels of the corresponding leaf vertices of the two trees are identical.

\begin{corollary}\label{tiplabeled}
$T_1, T_2\in\T_r^{\ell}$ are isomorphic if and only if $P_{\ell}(T_1)=P_{\ell}(T_2)$. 
\end{corollary}

\proof To prove this corollary, we claim that if a polynomial $P$ in $\Z[x,y]$ is an irreducible in $\Z[x,y]$ then the polynomial $Q$ in $\Z[x_1,x_2,...,x_t,y]$ by changing each $x$ in $P$ to some $x_i$ is also an irreducible in $\Z[x_1,x_2,...,x_t,y]$. Then, the corollary follows the proof of Theorem \ref{main}.  The proof of the claim is trivial because if a polynomial $Q\in\Z[x_1,x_2,...,x_t,y]$ is not an irreducible, say $Q=\prod_{i=1}^kQ_i$, then by substituting any $x_i$ with $x$ in the equation $Q=\prod_{i=1}^kQ_i$, we have  $P=\prod_{i=1}^kP_i$ where $P$ and $P_i$ are in $\Z[x,y]$ for all $1\leq i\leq k$. This contradicts that $P$ is an irreducible in $\Z[x,y]$. Hence, the polynomial $Q$ obtained by substituting $x$ in $P$ with some $x_i$ is an irreducible in  $\Z[x_1,x_2,...,x_t,y]$.\qed 

\bigskip
\noindent
Let $\T_u^{\ell}$ be the set of leaf labeled unrooted trees and define $\T^{\ell}=\T_u^{\ell}\cup \T_r^{\ell}$. Since the polynomial $P_{\ell}:\T_r^{\ell}\to \Z[x_1,x_2,...,x_t,y]$ is a complete isomorphism invariant for leaf labeled rooted trees and for any tree $T$ in $\T_r^{\ell}$, $P_{\ell}(T)$ is an irreducible in the polynomial ring, according to Section \ref{r4}, we can generalize the polynomial to a polynomial $P_{\ell}:\T^{\ell}\to \Z[x_1,x_2,...,x_t,y]$ such that for any leaf labeled rooted tree, its polynomial is defined as in Definition \ref{polyll} and for any leaf labeled unrooted tree $T$, $P_{\ell}(T)=\prod_{T_i\in\R_T}P_{\ell}(T_i)$. 

\begin{corollary}\label{tiplabeled2}
$T_1, T_2\in\T^{\ell}$ are isomorphic if and only if $P_{\ell}(T_1)=P_{\ell}(T_2)$. 
\end{corollary}
\proof To prove this corollary, we only need to generalize Lemma \ref{ll} for leaf labeled trees, that is, $T, B\in\T_u^{\ell}$ are isomorphic if and only if there exists a bijection $h:\R_{T}\to\R_{B}$ such that for any $T_i$ in $\R_T$, $T_i$ is isomorphic to $h(T_i)$. Note that if $\phi:T\to B$ is the isomorphism, then any leaf edge $e$ of $T$ should have the same label as the leaf edge $\phi(e)$ of $B$.  Besides, for any $T_i\in\R_T$, if $T_i$ is constructed by contracting a leaf edge $e$ of $T$, we consider that the root vertex of $T_i$ is labeled and it is of the same label as the leaf edge $e$ of $T$. Moreover,  $T_i\simeq h(T_i)$ requires not only the corresponding leaf vertices but also the root vertices to have the same label. Thus, this can be proved similarly to the proof of Lemma \ref{ll}. \qed

\section*{Acknowledgments} 

\noindent
The author would like to thank Priscila Do Nascimento Biller, Caroline Colijn and G\'abor Hetyei for helpful comments. This work was supported by the grant of the Federal Government of Canada's Canada 150 Research Chair program to Dr.~Caroline Colijn.

\bibliographystyle{plain}

\end{document}